\magnification=\magstep1
\input amstex
\documentstyle {amsppt}

\def\section#1{\bigskip\goodbreak\noindent{\bf #1}}

\long\def\Comment#1\endComment{\relax}

\def\about{\sim}

\def\rk{\operatorname{rank}}

\def\Z{\Bbb Z}

\def\set{\ S}

\NoBlackBoxes
\def\1{\text{1\!\!1}}
\document
\nopagenumbers 
\vskip1cm
\topmatter
\title
 Upper Bounds for the Davenport Constant
\endtitle
\author
R. Balasubramanian
and Gautami Bhowmik
\endauthor
\abstract
We prove that for all but a certain number of abelian groups of order $n$
the Davenport constant is atmost $\frac n k+k-1$ for positive integers $k\le 7$.
 For groups of rank three  we improve on the existing bound involving the
Alon-Dubiner constant.

\endabstract

\endtopmatter
\vskip-.25in

\document
\footline{\eightpoint The \the\month.\the\day.\the\year\hfill}

\section{I. Introduction}

Let $G$ be an abelian group of order $n$. A sequence  of elements (not necessarily
distinct) of $G$ is called a zero sum sequence of $G$ if the sum of its
components  is 0. The  zero-sum constant $ZS(G)$ of $G$ is defined to be the
smallest integer $t$ such that every sequence of length $t$ of $G$ contains a
zero-sum subsequence of length $n$, while  the Davenport constant $D(G)$ is  the
smallest integer $d$ such that every sequence of length $d$ of $G$ contains a
zero-sum subsequence.

The study of the zero-sum constant dates back to the Erd\" os-Ginzburg-Ziv
theorem of 1961 [EGZ]. On the other hand 
Davenport in 1966 introduced $D(G)$ as the maximum possible number of prime ideals
(with multiplicity) in the prime ideal decomposition of an irreducible element
of the ring of integers of an algebraic number field whose ideal class group
is $G$. More recently, Gao [G] proved that these two constants are closely
related,
i.e. $ZS(G)=|G| +D(G)-1$. It is thus enough to study any one of these
constants.

Apart from their interest in zero sum problems of additive number theory 
and non-unique factorisations in algebraic number theory, these constants play
an important role in
graph theory (see, eg, [Ch]). However their determination is still an open problem.

We consider the cyclic decomposition of a group of  $\rk r$,  i.e. 
 $G\sim {\Z}_{d_1}\oplus {\Z}_{d_2}\oplus\cdots
\oplus {\Z}_{d_r}$, where $d_i$ divides $d_{i+1}$ . It is clear that 
$M(G)=1+\sum_{i=1}^r  (d_i-1)$ is a lower bound for $D(G)$.

It was proved that $D(G)=M(G)$ for $p$ groups and for groups of rank 1 or 2,
independently by Olson [O] and Kruswijk [B1] and the equality is also known to
hold for several other groups. Olson and Baayen both conjectured that the
equality holds for all finite abelian groups. The conjecture however turned out to be
false . Geroldinger and Schneider [GS] in 1992 in fact showed that
for all groups of rank greater than 3, there exist infinitely many cases
where $D(G)>M(G)$. 

As far as upper bounds are concerned, the Erd\" os-Ginzburg-Ziv theorem that
asserts that for a finite abelian group of order $n$,
$ZS(G)\le 2n-1$ [EGZ] has been improved. Alon, Bialostocki and
Caro [cited in OQ] proved that  $ZS(G)\le 3n/2$ if $G$ is non-cyclic.
Caro improved this bound to $ZS(G)\le 4n/3 +1$ if $G$ is neither cyclic nor
of the form ${\Z}_{2}\oplus {\Z}_{2t}$. On excluding ${\Z}_{3}\oplus {\Z}_{3t}$
as well, Ordaz and Quiroz [OQ] tightened the bound to $5n/4+2$. It is easy to
see that though it is true for  $k=1,2,3$  and $4$; for a general positive integer
$k$ we cannot say that
 $D(G)\le \frac {n}{k}+(k-1)$ 
whenever $G$ is not of the form ${\Bbb Z}_{u}\oplus {\Bbb Z}_{ut}, u< k$.

On the other hand , Alford, Granville and Pomerance [AGP] in 1994 used the bound 
$D(G)\le  m(1+\log \frac n{m}) $, where $m$ is the exponent of $G$,
 to prove the existence of infinitely many Carmichael numbers.

In this paper, we combine the two types of upper bounds to prove that
\bigskip 
\noindent{\it{Theorem .}} If $G$ is an abelian group of order $n$ and exponent
$m$, then for $k\le 7$, its  Davenport
constant $D(G) \le \frac {n}{k}+(k-1)$ whenever $\frac n{m}\ge k$ .
\bigskip
Thus when the ratio $\frac n{m}$ is small, we get an improvement on the [AGP] bound.

We expect the above result to be true for all $k\le\sqrt n$.
\bigskip
To study the Davenport constant, it is sometimes useful to use another constant
$D^{s}(G)$ which is the smallest integer $t$ such that every sequence of $G$
with length $t$ contains a zero sum subsequence of length atmost $s$.

Olson calculated  $D^p ({\Z}_{p}\oplus {\Z}_{p})$ for a prime number
$p$ and used it to
determine
$D(G)$ for the rank 2 case. As yet,  no precise result is known for  
$D^p ({\Z}_{p}^r)$ for $r\ge 3$. But Alon and Dubiner [AD] proved a remarkable
bound in 1995,
i.e. $D^p ({\Z}_{p}^r)\le c(r)p$. In fact  $c(r)$ can be taken to be
$(cr\log r)^r$where $c$ is an absolute constant. Dimitrov [D] used  the
Alon Dubiner constant to prove that $D(G)\le M(G)(Kr\log r)^r$ for an absolute
constant $K$.
In the general case we have only a slight improvement of Dimitrov's result. It is 
for the rank 3 case that our result is interesting.
\bigskip
\noindent{\it Theorem .} If  $G\about {\Bbb Z}_{a_1}\oplus {\Z}_{a_1a_2}\oplus
{\Z}_{a_1a_2a_3}$, we have
$$D(G)\le M(G)(1+\frac{K}{a_2a_3}),$$
where $K$ is a constant of the same order of magnitude as that obtained by
Alon-Dubiner.
\bigskip
At the end we give an elementary proof of a result of Alon Dubiner that 
helped them obtain the bound for $D^p(\Z_p^r)$.

\section{II. A General Bound}

We first prove a lemma which would help us find bounds for the Davenport constant when reasonable bounds
can be found for $D^s(G)$ and when $D(\Z_s ^3)$ can be calculated, for example when $s$ is a power of a prime.

\bigskip
\noindent{\it Lemma 1.} Let $D^s(\Z_s ^3)\le A$ , $u=[\frac {A-s}{D(\Z_s ^3)}]$, and let
$$
\matrix &h=&h_{a,b}&=D({\Z}_{a}\oplus {\Z}_{ab}),\qquad &a\ne 1 \\
&\ & &=D({\Z}_{b}), \qquad &a=1.
\endmatrix
$$ 
Then, if $h\ge u+1$, 
$$D(\Z_s\oplus{\Z}_{sa}\oplus {\Z}_{sab})\le B(h_{a,b}),$$
where  $B(h_{a,b})=(h_{a,b}-u-1)s+A.$ 
\bigskip

\newpage

{\it{Proof.}}
\bigskip
Let $\set $ be a set of $B(h)$ elements of $\Z_s\oplus{\Z}_{sa}\oplus
 {\Z}_{sab}$
 . Every sequence of length atleast $D^s({\Z}_s\oplus{\Z}_{s}\oplus
 {\Z}_{s}) $
contains a zero sum subsequence of length atmost $s$. Thus $B(h)$ contains one 
 zero sum subsequence of length atmost $s$. On removing this zero sum
 sequence, we would still have more than  $D^s({\Z}_{s}^3)$ elements left
 in $B(h)$.  Thus  there exist
disjoint subsets $A_1, A_2,\cdots A_{h-u-1}  $ in $S$ such that $|A_j|\le s$ and 
the sum
of the elements of $A_j$ is $(0,0,0)$ in $\Z_{s}^3.$  If these sets are removed from $B(h)$,
we still have more than $B(h)-(h-u-1)s\ge D^s(\Z_{s}^3)$ elements 
from which we can extract another subset $A_{h-u} $ disjoint from the others
 of length $\le s$ and still of
sum  $(0,0,0)$ in $D^s({\Z}_{s}^3)$. Now

$$B(h)-(h-u)s\ge A-s\ge uD({\Z}_{p}^3).$$
 So we can extract $u$ more subsets $A_{h-u+1},\cdots, A_h $
disjoint from the rest the sum of whose elements is still zero in ${\Z}_{s}^3$.

Thus we have $h$ disjoint subsets whose sum in ${\Z}_{s}^3$ is $(0,0,0)$,
i.e. the sum is of the form  $(a_j s, b_j s, c_j s)$ .
 Suppose that $a\ne 1$ and for $j\le h$ let 
$C_j=( b_j , c_j )$. Now $a_js$ is $0$ in  ${\Bbb Z}_{s}$
and since we have taken the sum over $h$ sets,
 there exists a
subcollection of $C_j$  whose sum is $(0,0)$ in 
${\Z}_{a}\oplus {\Z}_{ab}.$ The corresponding 
subcollection of $A_j$ will suit our purpose in ${\Z}_s\oplus{\ Z}_{sa}\oplus {\ Z}_{sab}$.

If $a=1$, we take $C_j=( c_j )$ and proceed as before.

\hfill$\square$
\bigskip
To get precise bounds it is often necessary to actually evaluate $D({\Z}_{s}^3)$
or atleast find reasonable bounds.
This is possible for small values of $s$ as follows :
\bigskip
\noindent{\it Lemma 2.}
We have,
$$\matrix D^s({\ Z}_s\oplus{\Z}_{s}\oplus {\Z}_{s})&=8,\  s=2\\ 
 &=17, \ s=3,\\&= 22,\ s=4.\endmatrix.$$

{\it{Proof.}} The first two assertions can be verified directly. We notice that any 9 
distinct elements in  ${\Z}_{3}^3$ contain a zero sum subsequence. The third
follows essentially from Harborth [H].

\Comment
For $s=4$, consider the 7 following elements : $x_1=(2,0,0), x_2=(0,2,0), x_3=(0,0,2),
x_4=(2,2,0), x_5=(2,0,2), x_6=(0,2,2)$ and $x_7=(2,2,2).  $

Let $A_i=\{y : 2y=x_i\}$. 
with $|A_i|=8$. Let $B_i=A_i\cup \{x_i\}$ and if possible let $S$ be a set of 43 elements with no zero sum
subsequence.

Now consider $C_i=S\cap B_i.$ Since $\cup_{i=1}^7B_i$ has all the non zero elements of
$\Z_4^3$, we have  $|\cup_{i=1}^7C_i|=43$. Thus there exists an $i$ such that $|C_i|\ge 7$ and $c_i$ has 
no zero sum subsequence.

We shall show that this leads to a contradiction. We take $i=1$ for convenience, the arguments 
would work for a general $i$. Thus $$A_1=\{(1,0,0), (3,0,0),(1,2,0),(3,2,0),(1,0,2), (3,0,2),(1,2,2), (3,2,2)\}.$$

Since $C_1$ is zero sum free, the following pairs of elements cannot occur in $C_1$ :
$((1,0,0), (3,0,0)), ((1,2,0),(3,2,0)), ((1,0,2), (3,0,2)), ((1,2,2), (3,2,2)).$

Thus $C_1$ has at most 5 distinct elements. Further the multiplicity of any element in $C_1$ is atmost 3
and $x_1$ can occur atmost once. Atmost one element can have multiplicity 2, since $2y_1+2y_2=2x_1=0.$
Further, if an element $y_1$ has multiplicity greater than 1, then $2y_1+x_1=0$ and $x_1$ does not belong to
$C_1$.

So we see that there exists no such $C_i$, and hence no such $S$.
\endComment
\hfill$\square$
\bigskip
Sometimes we cannot find an effective bound for  $D({\Z}_{s}^3)$ but we might be able to use the
following weaker bound which can be proved in the same way as Lemma 1.
\bigskip
\noindent{\it Lemma 3.} We have 
$$
D(\Z_{s^a}^{r-1}\oplus \Z_{s^a t})\le D(\Z_{s^a}^{r})t.$$
\bigskip
\Comment
For some estimates, it is necessary to use the constant $D_k(G)$ which is the smallest integer $t$
such that every sequence of $t$ elemnts of $G$ contains $k$ disjoint zero sum sequences.
\bigskip
\noindent{\it Lemma 4.} We have $D_2(\Z_3^3)\le 13$.
\bigskip
{\it Proof.} Consider the set $S=\{x_1,\cdots,x_{13}\}.$ Now either $S$ has a zero sum sequence of length less
than 7, in which case it has another such sequence among the remaining $D(\Z_3^3)$ elements. Or $S$ contains
no zero sum sequence of length upto 6. Consider the set  $T_i=\{x_1,\cdots,x_6,x_i\}$ with $ 7\le i \le 13.$
Now $T_i$ contains a zero sum sequence of length exactly 7. Thus 
$\sum_{j=1}^6 x_j+x_i=0$ where  $7\le i \le 13.$ This gives $x_7=\cdots=
x_{13}$ and $x_7+x_8+x_9=0$ in $\Z_3^3.$

\hfill$\square$
\bigskip
\noindent{\it Remark.} By the same method we could prove that $D_k(\Z_p^a)\le k(D(\Z_p^a)-1)+1.$
\endComment
\bigskip
\noindent{\it Theorem 1.} If $G$ is an abelian group of order $n$ and exponent
$m$, then for every positive integer $k\le 7$, its  Davenport
constant $D(G)$ is atmost $\frac n{k}+(k-1)$ whenever $\frac n {m}\ge k$ .
\bigskip
{\it{Proof.}} We notice that the exceptions to the bounds
stated in the theorems of  Erd\" os-Ginzburg-Ziv [EGZ], Alon-Bialostocki-Caro
[ABC], 
Caro [C] and  Ordaz-Quiroz [OZ]
 can be reformulated as the cases where  $\frac n {m}\ge k$ to assert 
our result for k= 1,2,3 and 4 
 respectively.  

It is known  [AGP] that
$$D(G)\le m(1+\log \frac n{m})$$
and the condition $ m(1+\log \frac n{m})\le \frac n{k}+k-1$ is satisfied whenever $\frac n {m}\ge 31$
for $ k=7$. Thus
it suffices to examine the groups where $\frac n {m}\le 31$ .
\smallskip
Case 1 :   $\rk (G)\ge 5$.

We notice that for a group of rank greater than 5, $\frac n {m}$ is always greater than 31.
 Let  $$G\sim {\Bbb Z}_{a_1}\oplus {\Bbb Z}_{a_1a_2}\oplus
{\Z}_{a_1a_2a_3}\oplus {\Bbb Z}_{a_1a_2a_3a_4}
\oplus {\Z}_{a_1a_2a_3a_4a_5}.$$

Here $n=a_1^5a_2^4a_3^3a_4^2a_5$, and $m=a_1a_2a_3a_4a_5$. 
Since $a_1\ge 2$, $\frac n m\le 31$ only when $a_1=2,a_2=a_3=a_4=1$. 
Now, a result of [OQ] says that for any abelian group $K$,  
$$D({\Bbb Z}_{2}\oplus {\Bbb Z}_{2}\oplus{\Bbb Z}_{2}\oplus K)\le 2D(K)+3.$$

Taking $K$ to be ${\Bbb Z}_{2}\oplus{\Bbb Z}_{2t}$, we get
$D(G)\le 4t+5\le \frac{ 32}{k}t+k-1$ for $k=5,6,7.$
\smallskip
Case 2 :  $\rk (G)=4$. 

The condition $\frac n {m}=a_1^3a_2^2a_3 < 31$ is satisfied only for the
following groups of rank 4 
that would violate the
AGP condition would be of the form 
$$G_1\sim {\Z}_{2}\oplus {\Z}_{2}\oplus{\Z}_{2}\oplus {\Bbb Z}_{2t},$$
$$G_2\sim {\Z}_{2}\oplus {\Z}_{2}\oplus{\Z}_{4}\oplus {\Z}_{4t},$$
$$ G_3\sim  {\Z}_{2}\oplus {\Z}_{2}\oplus{\Z}_{6}\oplus {\Z}_{6t},$$
and
$$G_4\sim {\Bbb Z}_{3}\oplus {\Bbb Z}_{3}\oplus {\Bbb Z}_{3}\oplus
{\Bbb Z}_{3t}.$$ 
 However the first case satisfies the
stronger condition of the Baayen-Olson conjecture, i.e. $ D(G)=M(G)$. This was
proved for odd $t$ [B2] and for even $t$ it follows from the fact that in this case
$$G_1=H\oplus{\Bbb Z}_{p^ku}$$
$H$ being a $p$-group and $p^k\ge M(H)$, a case that satisfies the BO
conjecture [vE]. 

For $G_2$ we split it as a sum of two groups $H$ and $K$ and use the estimate
(see eg [C]),
$$D(H+K)\le (D(H)-1)|K| + D(K). $$

We take $H$ to be ${\Z}_{2}\oplus 
{\Z}_{4}\oplus {\Z}_{4t}$. Then $D(H)=M(H)$ (see [vE]). Thus 
$D(G_2)\le 8t+8$ which is  less than $\frac n k+k-1$ for all $t$ when $k=5$ and
for $t>1$ when $k=6,7$. But for $t=1$ we have a $p$-group.

The same argument works for $G_3$. For $G_4$ we use Lemma 3 and get
$$D(G_4)\le 9t\le \frac{81}{7}t+6$$ for $k=7$. Since $\frac n m=27$ the
inequality is already satisfied by the AGP bound for $k=5,6$. 
\smallskip
Case 3 :  $\rk (G)=3.$ 

Since $a_1^2a_2\ge 31$ ensures that 
$D(G)\le \frac {n}{k}+k-1 $, and $\frac n {m} \ge k$ we are left with the cases 
$G_5\sim {\Bbb Z}_{2}\oplus {\Bbb Z}_{2u}\oplus
{\Bbb Z}_{2ut},\ 1< u< 8$  
 , $G_6\sim {\Bbb Z}_{3}\oplus {\Bbb Z}_{3v}\oplus
{\Bbb Z}_{3vt}, v=1,2,3$ ;  $G_7\sim {\Bbb Z}_{4}\oplus {\Bbb Z}_{4}\oplus
{\Bbb Z}_{4t}$ and  $G_8\sim {\Bbb Z}_{5}\oplus {\Bbb Z}_{5}\oplus
{\Bbb Z}_{5t}$.

Now  $G_5$  satisfies the BO conjecture.  This follows from the fact
 that $u$ has no prime divisor greater than 11, which is a sufficient
condition from a result 
  of van Emde Boas [vE].

With $s=3,a=1, b=t$ in Lemmas 1 and 2, we obtain, for $k=5$ or $6$ that 
$$D({\Bbb Z}_{3}\oplus {\Bbb Z}_{3}\oplus
{\Bbb Z}_{3t})\le 3t+8\le \frac{27}{k}t+k-1,$$ whenever $t\ge 2$.

When $t=1$, we have a $p$-group and the BO conjecture is satisfied. 

For $k=7$ we know that for ${\Bbb Z}_{3}\oplus {\Bbb Z}_{3}\oplus
{\Bbb Z}_{6}$  the BO condition is realised [vEK]
and we are within the claimed bound. The same is true for 
the cases $v=2,3 $ in $G_6$ [vE].
\Comment
We consider a set of 13 elements of $\Z_3^2\oplus\Z_6$.  From Lemma 4 we know that it
 contains two disjoint zero sum sequences in $\Z_3^3$, say $T_1$ and $T_2$.
Let their sums be $(3a,3b,3c)$ and $(3k,3l,3m)$ respectively. Now if either $c$ or $m$ is even,
we obtain a zero sum sequence in $\Z_3^2\oplus\Z_6$, while if both are odd, $T_1\cup T_2$ is the required
zero sum sequence.

For the case $v=2,3 $ in $G_6$, the BO condition is realised [vE]
and we are within the claimed bound.
\endComment
For $G_7$ we use Lemmas 1 and 2 with $s=4, a=1,b=t$ to obtain that

$$D(G_7)\le 4t+27\le \frac{64}{7}t+6$$ for $t>4, k=7$.  Lemma 3  gives the desired
bound for $k=5$ or 6 ,  $t\le 3, k=7$ in $G_7$ as well as for all cases of  $G_8$.
\medskip
Case 4 :  $\rk (G)=2$. 

It is well known that  $D(G)=a_1+a_1a_2-1$ and
the inequation $$a_1+a_1a_2-1\le \frac {a_1^2a_2}{k}+k-1$$ is always true for $a_1=\frac n {m}\ge k$.

\hfill$\square$
\bigskip
\noindent{\it Remark. }  This bound is tight, since 
$D({\Bbb Z}_{k}\oplus {\Bbb Z}_{kt})=kt+k-1$.

\bigskip
\noindent{\it Conjecture. }  We believe that Theorem 1 is true for all $k\le\sqrt n$. Notice that this is a 
weaker claim than the Narkiewicz-\'Sliwa conjecture that $D(G)\le M(G)+r-1$ for a group of rank $r$.

\section{III. Rank 3 case} 

We now use the Alon-Dubiner theorem for improving the existing bound 
for the
Davenport constant when the rank of the group is 3 which is [D] 
$$D(G)\le K_3M(G),$$
where $K_3$ is a constant of the same order of magnitude as
that of Alon-Dubiner, and $M(G)=a_1a_2a_3 +a_1a_2 +a_1-2$. Our method also gives a minor
improvement
for the higher rank cases.

We state a Lemma which can be seen as a generalisation of Olson's result for
the rank 2 case.
 
\bigskip
\noindent{\it Lemma 5.} 
 Let $d$ be a  divisor of $a$ and let 

$$ h=D({\Z}_{a/d}\oplus {\Z}_{ab/d}\oplus {\Z}_{abc/d}), a\ne
d,$$
$$ =D({\Z}_{b}\oplus {\Z}_{bc})\ ,\qquad a=p\ ,\  b\ne 1,$$
$$ =D({\Z}_{c})\ ,\qquad a=d, b=1, c\ne 1.$$
Then 
$$D(\Z_a\oplus{\Z}_{ab}\oplus {\Z}_{abc})\le B(h),$$
where  $B(h)=(h-u-1)d+A,$ and $A$ and $u$ are as defined in Lemma 1.
\bigskip
\noindent{\it{Proof.}} Same as that of Lemma 1.
\hfill$\square$
\bigskip
\noindent{\it Theorem 2.}
 Let $G\sim {\Bbb Z}_{a_1}\oplus {\Bbb Z}_{a_1a_2}\oplus
{\Bbb Z}_{a_1a_2a_3}$. Then
$$D(G)\le a_1a_2a_3 +a_1a_2 +Ka_1,$$
where $K$ is a constant of the same order of magnitude as that of Alon-Dubiner.
\bigskip

{\it Proof . }

We use Lemma 3 above and the Alon Dubiner bound,
$$D^p({\Z_p}^r)\le c(r)p,$$
where $c(r)$ is a constant.
In particular, for $r=3$, we write $D^p({\Z_p}^3)\le (K+3)p$ with
$(K+3)p\ge 7p-4$. Thus $hp+Kp\ge hp+4p-4$.

For fixed $a_2,a_3$  we write 
$h(a_1)=D({\Bbb Z}_{a_1}\oplus {\Bbb Z}_{a_1a_2}\oplus{\Bbb Z}_{a_1a_2a_3}).$

Using Lemma 3 we see that if $p$ divides $a_1$ , 
$$ h(a_1)\le h((a_1/p)+K)p.$$
 Let $a_1=p_1p_2\cdots p_t$ with
$p_i\ge p_{i+1}$.
Thus $$ h(p_1p_2\cdots p_t)\le h((p_2\cdots p_t)+K)p.$$
Repeating the above process we get
$$h(a_1)\le a_1h(1)+K(p_1p_2\cdots p_t+p_1p_2\cdots
 p_{t-1}+\cdots+p_1)
$$
But $$p_1p_2\cdots p_t+p_1p_2\cdots
 p_{t-1}+\cdots+p_1 =a_1(1+\frac 1{p_t}+\frac 1{p_tp_{t-1}}+\cdots
+\frac 1{p_tp_{t-1}\cdots p_2})\le 2a_1.$$

This gives $D({\Z}_{a_1}\oplus {\Z}_{a_1a_2}\oplus{\Z}_{a_1a_2a_3})
\le a_1D( {\Z}_{a_2}\oplus{\Z}_{a_2a_3}) +2Ka_1,$

i.e. $$D({\Z}_{a_1}\oplus {\Z}_{a_1a_2}\oplus{\Z}_{a_1a_2a_3})\le a_1a_2a_3 +a_1a_2 +(2K-1)a_1.$$

\hfill$\square$
\bigskip
{\it Remark.} For the case of a general $r$ we get $D(G)\le M(G)(1+\frac{K_r}{a_{r-1}a_r})$
and the improvement from the existing bound comes into picture only when 
$a_{r-1}$ and $a_r$ are large.

The proof of Theorem 2 uses an inequality of [Proposition 2.4,AD]. Here we give a slightly
improved constant for the inequality. The proof goes on the same lines as
[AD] but uses no graph theory. We include it here for the sake of completion.
\medskip
{\it Theorem 3.} Let $A$ be a subset of ${\Bbb Z_p^d}$ such that no
hyperplane contains more than $\mid A \mid /4W$ elements of $A$. Then for all
subsets $Y$ of ${\Bbb Z_p^d}$ containing at most $p^d/2$ elements, there is an 
element $a\in A$ such that
$$\mid (a+Y)\backslash Y | \ge \frac W{5p} \mid Y \mid. $$ 
\bigskip
{\it Proof.} If possible, let there exist no such $a$. Then 
$L(a)=|(a+Y)\backslash Y | \le \frac W {5p} |Y|  $ for all $a\in A$. Since $L(ja)\le jL(a)$, we get
 $L(ja)\le \frac {jW}{5p} \mid Y \mid $ for all $j\le p/W$.

This gives  $$M(ja)=L(ja)+L(-ja)\le \frac  {2jW}{5p} \mid Y \mid. $$ 

Let $J=[\frac p W]$. 
Then
$$S=\sum_a\sum_{1\le j \le J}M(ja)\le J(J+1) \frac W{5p}|Y||A|.$$

On the other hand we shall get a lower bound for $S$. For any $b$ define
$$T(b)=\frac{1}{|G|}\sum_{x}(1-l^{\bar b.\bar x})
|\sum_y l^{\bar x .\bar y}|^2,$$
where for notational convenience we write $ l$ for  $e^{\frac{2i\pi}{p}}$ and $G$ for the group $\Z_p^d$.

Then$$T(b)=\frac{1}{|G|}\sum_{x}(1-l^{\bar b.\bar x})
\sum_{y_1,y_2} l^{\bar x .(\bar y_1-\bar y_2)}.$$

$$=\frac{1}{|G|}\sum_{y_1,y_2}(\sum_{x}l^{\bar x .(\bar y_1-\bar y_2)}-\sum_{x}l^{\bar x .(\bar y_1-\bar y_2-\bar b)})
 .$$
$$=B-D.$$

Clearly $B=|Y|$ and $D$ is the number of solutions of the equation $\bar y_1-\bar y_2=\bar b$
which is the same as $(b+Y)\cap Y.$

Thus $B-D=|(b+Y)\backslash Y|=L(b)$.
Thus  $$M(ja)=L(ja)+L(-ja)=T(ja)+T(-ja)=
\sum_{x}(2-l^{j\bar a.\bar x}-l^{-j\bar a.\bar x})
|\sum_y l^{\bar x .\bar y}|^2.$$
$$=\frac{4}{|G|}\sum_{x}\sin ^2(\frac{\pi }{p}j\bar a.\bar x)|\sum_y l^{\bar x .\bar y}|^2.$$
Then 
$$ S=\frac{4}{|G|}\sum_{x\ne 0}\sum_a\sum_j\sin ^2(\frac{\pi }{p}j\bar a.\bar x)
|\sum_y l^{\bar x .\bar y}|^2$$
$$\ge \frac{4}{|G|}\sum_{x\ne 0}R|\sum_y l^{\bar x .\bar y}|^2$$

where $R$ is a minorisation of $\sum_a\sum_j\sin ^2(\frac{\pi }{p}j\bar a.\bar x)$ for $x\ne 0$.

We then have $$
S\ge \frac{4R}{|G|}\sum _{x\in G}|\sum_y  l^{\bar x .\bar y}|^2-
\frac{4R}{|G|}\sum _{x=0}|\sum_y  l^{\bar x .\bar y}|^2 $$   
$$ \ge 4R|Y|-\frac{4R}{|G|}|Y|^2\ge2R|Y|,  $$
since $|Y|\le \frac{|G|}{2}$. On the other hand, to get a lower bound for $R$ we note that the 
least value is obtained by taking $j\bar a.\bar x$ as small as possible. Thus the condition that no hyperplane contains
more than $\frac {|A|}{4W} $ elements implies that $\bar a.\bar x\ge W$ for atleast $\frac{|A|}{2}$ values of $a$.
Considering only these values, we have $R\ge\frac{|A||J|}{8}$ and $S\ge\frac{|A||J|}{4}|Y|$.
This gives a contradiction.
\bigskip
\noindent{\it Acknowledgement.}  The authors are grateful to CEFIPRA (Project 2801-1)
for their financial support  for visits to each others' institute. 
\bigskip
\centerline{{\rm REFERENCES}}

\def\ptt{\hskip 0pt plus 2.5pt minus 0.5pt}
\long\def\Comment#1\endComment{\relax}
\def\Ref#1\endref{%
\def\nobox{}%
\def\paperbox{}%
\def\bybox{}%
\def\jourbox{}%
\def\bookbox{}%
\def\yrbox{}%
\def\pagesbox{}%
\def\volbox{}%
\def\otherbox{}%
\def\submittedbox{}%
\def\no##1{\def\nobox{\hbox to 23pt{\rm ##1\hfill}}}%
\def\paper##1{\def\paperbox{\ptt\sl ``##1''\ptt\ }}%
\def\by##1{\def\bybox{\ptt\sl\  \rm ##1\ptt\ }}%
\def\jour##1{\def\jourbox{\ptt\sl\   \rm ##1 \ptt}}%
\def\book##1{\def\bookbox{\ptt\sl\ publi\'e par \rm ##1 \ptt}}%
\def\other##1{\def\otherbox{\ignorespaces~,\ptt\ \rm ##1 \ptt}}%
\def\submitted##1{\def\submittedbox{\ptt\sl\ soumis \`a \rm ##1 \ptt}}%
\def\yr##1{\def\yrbox{\ptt\ (\bf ##1)}}%
\def\pages##1{\def\pagesbox{~,\ptt\sl\ pages ##1}}%
\def\vol##1{\def\volbox{\ptt\rm ##1\ptt}}%
#1%
\noindent\nobox%
\paperbox%
\bybox%
\jourbox%
\bookbox%
\submittedbox%
\otherbox%
\volbox%
\yrbox%
\pagesbox~.%
}
{\eightpoint
\Ref
\no{[ABC]}
\paper {Extremal zero sum problem}
\by {}
\by {N.~Alon, A.~Bialostocki \& Y .~Caro}
\jour {Manuscript,\goodbreak cited in [OQ]}
\vol {}
\endref

\Ref
\no{[AD]}
\paper {A Lattice Point Problem and Additive Number Theory }
\by {}
\by {N.~Alon \& M.~Dubiner}
\jour {Combinatorica}
\vol {15}
\yr{1995}
\pages {301-309}
\endref

\Ref
\no{[AGP]}
\paper {There are infinitely many Carmichael numbers}
\by {}
\by {W.R.~Alford, A.~Granville \& C .~Pomerance}
\jour {Annals of Math}
\vol {}
\yr{1994}
\pages {703-722}
\endref

\Ref
\no{[B1]}
\paper {Een combinatorisch probleem voor eindige Abelse groepen}
\by {}
\by {P.C~ Baayen}
\jour {Colloq. Discrete Wiskunde Caput 3, Math Centre, Amsterdam}
\vol {}
\yr{1968}
\endref

\Ref
\no{[B2]}
\paper {$(C_2 \oplus C_2 \oplus C_2 \oplus C_{2n})!$ Is True for Odd n}
\by {}
\by {P.C~ Baayen}
\jour {Report ZW-1969-006, Math Centre, Amsterdam}
\vol {}
\yr{1969}
\endref

\Ref
\no{[C]}
\paper {On zero sum subsequences in abelian non-cyclic groups}
\by {}
\by { Y .~Caro}
\jour {Israel Jour. of Math}
\vol {92}
\yr{1995}
\pages {221-233}
\endref

\Ref
\no{[Ch]}
\paper {On the Davenport Constant, the Cross Number, and their Application in
Factorization Theory}
\by {in Zero-dimensional commutative rings}
\by {S-T.~Chapman }
\jour {Lecture Notes in Pure and Applied Maths,Dekker, eds Anderson \& Dobbs}
\vol {171}
\yr{1995}
\pages {167-190}
\endref

\Ref
\no{[D]}
\paper {Zero-sum problems in finite groups}
\by {}
\by {V.~Dimitrov}  
\jour {}
\vol {}
\yr{2003}
\endref

\Ref
\no{[EGZ]}
\paper {Theorem in Additive Number Theory}
\by {}
\by {P.~Erd\"os,A.~Ginzburg \& A.~Ziv}
\jour {Bull. Research Council Israel}
\vol {10}
\yr{1961}
\pages {41-43}
\endref

\Ref
\no{[G]}
\paper {Addition theorems for finite abelian groups}
\by {}
\by {W.~ Gao}
\jour {J. Number Theory}
\vol {53}
\yr{1995}
\pages {241-246}
\endref

\Ref
\no{[GS]}
\paper {On Davenport's Constant}
\by {}
\by {A.~ Geroldinger \& R.~Schneider}
\jour {J. Combinatorial Theory, Series A}
\vol {61}
\yr{1992}
\pages {147-152}
\endref

\Ref
\no{[H]}
\paper {Ein extremalproblem f\"ur Gitterpunkte}
\by {H. Harborth}
\jour {J. Reine Angew.Math.}
\vol {262}
\yr{1973}
\pages {356-360}
\endref

\Ref
\no{[O]}
\paper {A combinatorial problem on finite Abelian groups I, II}
\by {}
\by {J.E.~ Olson}
\jour {J. Number Theory}
\vol {1}
\yr{1969}
\pages {8-11, 195-199}
\endref

\Ref
\no{[OQ]}
\paper {The Erd\"os-Ginzburg-Ziv Theorem in Abelian non-Cyclic Groups}
\by {}
\by {O.~Ordaz \& D .~Quiroz}
\jour {Divulgaciones Matematicas}
\vol {8, no. 2}
\yr{2000}
\pages {113-119}
\endref

\Ref
\no{[vE]}
\paper {A combinatorial problem on finite Abelian groups II}
\by {}
\by {P.~van Emde Boas }
\jour {Report ZW-1969-007 Math Centre, Amsterdam}
\vol {}
\yr{1969}
\endref

\Ref
\no{[vEK]}
\paper {A combinatorial problem on finite Abelian groups III}
\by {}
\by {P.~van Emde Boas \& D.Kruyswijk}
\jour {Report ZW-1969-008 Math Centre, Amsterdam}
\vol {}
\yr{1969}
\endref
\bigskip
Authors' Address:

1.Institute of Mathematical Sciences, Chennai,India.Email:balu\@imsc.res.in

2.Universit\'e de Lille 1, 
   Laboratoire de Math. U.M.R. CNRS 8524,
  59655 Villeneuve d'Ascq Cedex, France.Email:bhowmik\@math.univ-lille1.fr
\end
\Ref
\no{[]}
\paper {}
\by {}
\by {.~ \& .~}
\jour {}
\vol {}
\yr{19}
\pages {}
\endref

\enddocument